\documentclass[letterpaper, 10 pt, conference]{ieeeconf}  

\IEEEoverridecommandlockouts                              

\usepackage{graphics} 
\usepackage{epsfig} 
\usepackage{mathptmx} 
\usepackage{times} 
\usepackage{amsmath} 
\usepackage{amssymb}  

\usepackage{subfig}

\usepackage{tikz}
\usepackage{tikzscale}
\usepackage{pgfplots}
\pgfplotsset{compat=1.12}
\usetikzlibrary{matrix}
\usepgfplotslibrary{groupplots}

\pgfplotsset{compat=newest}

\DeclareMathOperator{\atan}{atan}
\DeclareMathOperator{\asin}{asin}

\title{\LARGE \bf
Guaranteeing Consistency in a Motion Planning and Control Architecture Using a Kinematic Bicycle Model
}

\author{Philip Polack$^{1}$, Florent Altch\'e$^{1,2}$, Brigitte d'Andr\'ea-Novel$^{1}$ and  Arnaud de La Fortelle$^{1}$
\thanks{This work was supported by the international Chair MINES ParisTech - Peugeot-Citro\"en - Safran - Valeo on ground vehicle automation}
\thanks{$^{1}$Philip Polack, Florent Altch\'e, Brigitte d'Andr\'ea-Novel and  Arnaud de La Fortelle are with Center of Robotics, Mines ParisTech, PSL Reasearch University, 60 boulevard Saint-Michel, 75006 Paris, France 
{\tt\small \{philip.polack, florent.altche, brigitte.dandrea-novel, arnaud.de\_la\_fortelle\} @mines-paristech.fr}} 
\thanks{$^{2}$Florent Altch\'e is also with \'Ecole des Ponts ParisTech, Cit\'e Descartes, 6-8 Av Blaise Pascal, 77455 Champs-sur-Marne, France.}%
}

\begin{document}

\maketitle
\thispagestyle{empty}
\pagestyle{empty}

\begin{abstract}

This paper proposes to combine a 10Hz motion planner based on a kinematic bicycle Model Predictive Control (MPC) and a 100Hz closed-loop Proportional-Integral-Derivative (PID) controller to cope with normal driving situations. Its novelty consists in ensuring the feasibility of the computed trajectory by the motion planner through a limitation of the steering angle depending on the speed. This ensures the validity of the kinematic bicycle model at any time. The architecture is tested on a high-fidelity simulation model on a challenging track with small curve radius, with and without surrounding obstacles. 

\end{abstract}

\section{INTRODUCTION}
\label{sec:introduction}

As research on autonomous vehicles is getting more and more mature, the question of consistency between its different layers, namely perception, localization, planning and control, is becoming crucial to ensure the safety of the vehicle at all time. An ill-designed vehicle architecture might be very critical for its safety, even though each layer is well-designed independently.

In this paper, we propose to focus on the problem of consistency between planning and control. Motion planning and control problems are two different but highly related problems. The first consists in computing a feasible trajectory (in terms of vehicle's dynamics) for the vehicle avoiding the surrounding obstacles such as other vehicles, pedestrians, or non moving objects. The second is acting on the actuators, \textit{i.e.} the gas pedal, brake pedal and steering wheel, in order to track the trajectory obtained by the motion planner, while ensuring the stability of the system and, if possible, a smooth drive. Therefore, the properties of the design of each layer are quite different. 

\begin{figure}[thpb]
	\centering
	\includegraphics[scale=0.5]{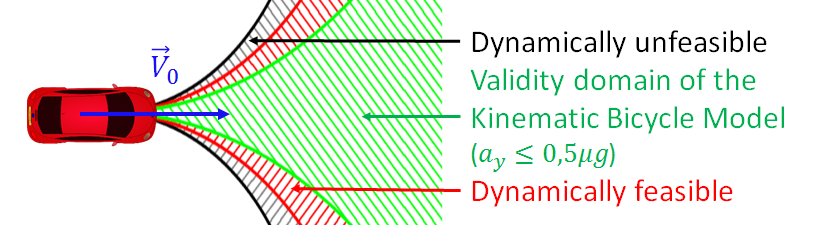}       
	\caption{Dynamic feasibility of motions considered as ``admissible" for a kinematic bicycle model at a given speed $V_0$ in the physical space.} 
	\label{fig:consistency_kbm}
\end{figure}

For motion planning purpose, the algorithm will have to explore the space of feasible solutions, which is computationally expensive. Thus, it requires both a rather simple model of the vehicle, such as a point-mass \cite{Nilsson2015} or a kinematic bicycle model \cite{Li2017}, and a rather low-frequency (around 5-10Hz). On the contrary, control is usually operating at high-frequency (around 100Hz) to ensure a good tracking of the reference trajectory and to be robust to modeling errors and disturbances. 
Moreover, the level of abstraction also differs: dealing with obstacles is one of the main task of the motion planner while the controller usually completely ignores them; the trajectory given by the motion planner is assumed to be safe within a certain margin and the goal of the controller is then to follow as well as possible the given trajectory, without considering the obstacles. Therefore, if the motion planner is planning outside the validity domain of its model, it might compute unfeasible trajectories from a dynamical point of view as shown in Figure~\ref{fig:consistency_kbm}. Then, the controller will not be able to follow the trajectories, putting the whole system in jeopardy! 

The problem of consistency, also referred to as ``proper modeling" in the literature \cite{Ersal2008}, arises with the need for simple models for motion planning: as simpler models are a priori less precise, there is a trade-off between the complexity (\textit{i.e.} the computational time) and the area of validity of the model. It has been studied in the literature by comparing simple models with more complex ones. For example, \cite{Liu2013} compares a motion planner and control MPC using a 2 Degrees of Freedom (2~DoF) linear-tire dynamic bicycle model with a 14~DoF vehicle model. The validity of this dynamic bicycle model is guaranteed by constraining the lateral acceleration to $0.5g$ on normal road conditions, condition which has been derived in \cite{Park2009}. However, the practical implementation of this system on a real vehicle is limited as it assumes a constant velocity. 
\cite{Kong2015} compares the performances of MPC using a kinematic bicycle model and one using a linear-tire dynamic bicycle model. The authors conclude that the kinematic bicycle model works better in most cases, except at high-speed, and suggests a further investigation on the impact of different lateral accelerations on the validity of the kinematic bicycle model. This has been done in our previous work \cite{Polack2017a} by comparing it to a high-fidelity 9~DoF vehicle model on constant radius curves. The conclusion was that the lateral acceleration should always remain lower than $0.5\mu g$, where $\mu$ is the road friction coefficient.

In the literature, many works can be found on motion planning and control architectures. They can be classified into two categories \cite{Berntorp2017}. In the first category, motion planning and low-level control are unified into one unique MPC formulation using different vehicle models, such as in \cite{Abbas2014}, \cite{HoujieJiang2016} and \cite{Liniger2015}. However, they are not robust to modeling errors, disturbances and parameter uncertainties as they lack a real-time feedback controller (the frequency of recomputing is rather low for low-level control). Thus, for stability reasons, it is preferable to use a high-frequency closed-loop controller. The second category separates motion planning from low-level control to overcome this problem (see for example \cite{Rosolia2017}, \cite{Li2017}). However, in that case, there are no guarantees on the feasibility of the planned trajectory as the model used for motion planning might be unable to capture the admissible motions of the vehicle. This problem has not been encountered yet as the number of tested kilometers with such architecture is small and the driving strategy chosen often very conservative, but it is very critical for the vehicle safety.

Therefore, this paper proposes a simple planning and control architecture for normal driving situations that guarantees the dynamic feasibility of the planned trajectory while being robust to parameter uncertainties and disturbances. The motion planner relies on a 10Hz kinematic bicycle Model Predictive Control approach, while the low-level closed-loop control signals are computed at a high-frequency (every 10ms) and are mainly based on Proportional-Integral-Derivative (PID) controllers. 
The motion planner is always planning a valid dynamic reference trajectory (green region of Figure~\ref{fig:consistency_kbm}) according to our previous work \cite{Polack2017a}, by constraining the lateral acceleration to be less than $0.5\mu g$. 
Our architecture is also less conservative than other ones: a local velocity planner based on the lateral acceleration constraint was designed in order to compute a heuristic speed for the motion planner, that pushes the vehicle's architecture towards its operational limits, while remaining safe. 


The paper is organized as following: first, the kinematic bicycle model with constrained lateral acceleration to ensure its validity is presented. Then, the consistent planning and control architecture chosen is explained in Section~\ref{sec:ctrl_archi}. Section~\ref{sec:SimModel} presents the 9~Degrees of Freedom (9~DoF) model used for our simulations. Section~\ref{sec:results} shows the results obtained in simulation with the proposed consistent planning and control architecture on a track with and without obstacles. Finally, Section~\ref{sec:CCL} concludes this paper.


\section{The kinematic bicycle model}
\label{sec:kbModel}

\begin{figure}[thpb]
	\centering
	\includegraphics[scale=0.55]{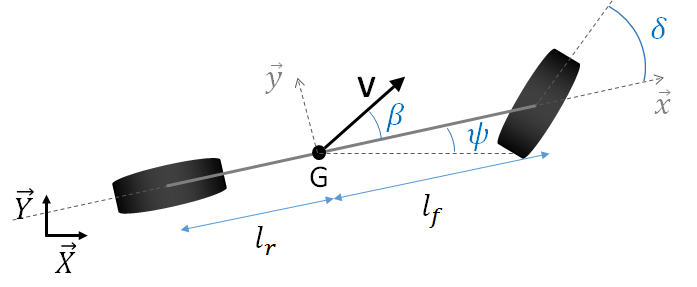}       
	\caption{Kinematic bicycle model of the vehicle.}
	\label{fig:kin_bic_model}
\end{figure}

In the kinematic bicycle model, the two front wheels (resp. the two rear wheels) of the vehicle are lumped into a unique wheel located at the center of the front axle (resp. of the rear axle) such as illustrated on Figure~\ref{fig:kin_bic_model}. In our case, the control inputs correspond to the acceleration $u_1$ of the vehicle and the steering rate $u_2$ of the front wheel. As the MPC computes piecewise constant inputs, this enables to have a smoother control because the reference speed and the steering angle become then continuous. Assuming that only the front wheel can be steered, the kinematic bicycle model\footnote{This is a misuse of language, as we are actually using the acceleration instead of the speed as input for better smoothness.} can then be written~\cite{Rajamani2012} as:
\begin{subequations}
	\label{eq:kin_bic_model}
	\begin{eqnarray}
	\dot{X} & = & V \cos \left(\psi + \beta(\delta) \right)\\ 		 \label{eq:kin_bic_model_1}
	\dot{Y} & = & V \sin \left(\psi + \beta(\delta) \right)\\ \label{eq:kin_bic_model_2}
	\dot{\psi} & = & \frac{V}{l_r} \sin\left(\beta(\delta)\right)\\\label{eq:kin_bic_model_4}
	\dot{V} & = & u_1\\ \label{eq:kin_bic_model_3}
	\dot{\delta} & = & u_2
	\label{eq:kin_bic_model_6} 
	\end{eqnarray}
	where $X$ and $Y$ are the global coordinates of the vehicle, $\psi$ the yaw angle, $V$ the speed and $\delta$ the front steering angle. The slip angle $ \beta$ at the center of gravity depends on $\delta$:
	\begin{eqnarray}
	\beta(\delta) & = & \atan\left(\tan(\delta) \frac{l_r}{l_f+l_r}\right) \label{eq:kin_bic_model_5}
	\end{eqnarray}
\end{subequations}

This model does not take into account the vehicle's dynamics such as slipping and skidding. Therefore, when used for planning problems, it might plan unfeasible trajectories (black region of Figure~\ref{fig:consistency_kbm}) which are not acceptable for the safety of the vehicle as explained in Section~\ref{sec:introduction}. However, in previous work \cite{Polack2017a}, we have shown that constraining the lateral acceleration of the vehicle under $0.5\mu g$ ensures the validity  of the model (green region of Figure~\ref{fig:consistency_kbm}).

\begin{figure}
	\centering
	\vspace{+0.1in}
	\includegraphics{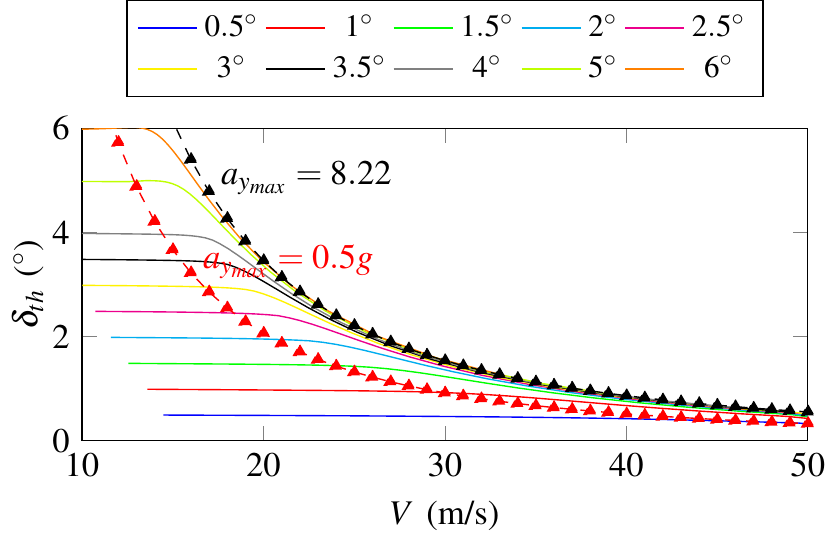}
	\caption{Comparison between the steering angle $\delta_a$ applied on a 9~DoF dynamic model of a vehicle (in color) and the one applied on a kinematic bicycle model $\delta_{th}$ in order to achieve a given turn radius $R$, at different speeds $V$ (cf. \cite{Polack2017a}).}
	\label{fig:ay_max}
\end{figure}

This constraint was obtained by comparing\footnote{The kinematic bicycle model remains valid as long as $\delta_a=\delta_{th}$ to achieve the same turn radius $R$, where $\delta_{th}$ is constant (does not depend on speed).} the steering angle $\delta_a$ applied on a complex 9~DoF dynamic vehicle model (presented in Section~\ref{sec:SimModel}) and the one applied on a kinematic bicycle model $\delta_{th}$, see Equation~(\ref{eq:delta_kin}), to follow a curve of radius $R$ at different speeds $V$. The results are displayed on Figure~\ref{fig:ay_max}. We chose curves as scenarios as they correspond to the situations where the vehicle is most likely to slip. 
\begin{eqnarray}
\label{eq:delta_kin}
\delta_{th} = \atan \left(\left(\frac{l_f}{l_r}+1\right)\tan \left(\asin \left(\frac{l_r}{R}\right)\right)\right)
\end{eqnarray}

\section{The planning and control architecture}
\label{sec:ctrl_archi}

\begin{figure}[thpb]
	\centering
	\includegraphics[scale=0.32]{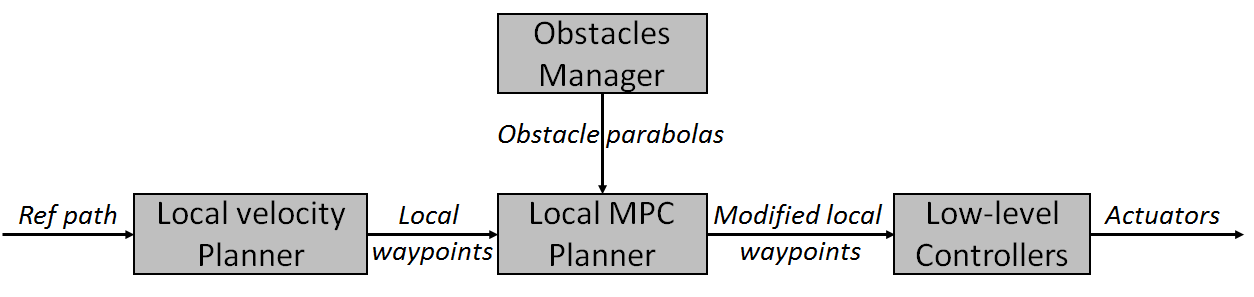}       
	\caption{Control architecture of the ego-vehicle.}
	\label{fig:control_archi}
\end{figure}

The goal of the planning and control architecture is to compute the actuator inputs (steering angle and wheel torques) to send to the vehicle, while setting a target speed and avoiding obstacles located on the path. Our system is composed of four main components to deal with most driving situations, namely a local velocity planner, an obstacle manager, a local MPC planner and low-level controllers (see Figure~\ref{fig:control_archi}). These parts are detailed respectively in Subsections~A-D. The main characteristics of our design are:
\begin{itemize}
	\item A smart computation by the local velocity planner of a heuristic speed $V_{heur}$ based on the lateral constraint for the kinematic bicycle model, to guide the MPC towards high-speed maneuvers.
	\item A combination of a 10Hz MPC planner for anticipating new situations, with a high-frequency (100Hz) closed-loop low-level controller for a better tracking of the reference trajectory.
	\item A validity guarantee of the kinematic bicycle model used in the MPC through the constraint on lateral acceleration, generating thus only feasible trajectories to follow for the low-level controllers. 
\end{itemize}

\subsection{The Local Velocity Planner}
\label{ssec:Vtarget}

The reference input of our architecture only contains the future coordinates $(x_i,y_i)_{i \in \mathbb{N}}$ of the center line of the road. However, we want to drive our vehicle at the highest possible speed that keeps our architecture safe. Therefore, a velocity term should be added to the cost function of the local MPC planner (see Section~\ref{ssec:planner}) that fosters higher but safe speeds. Choosing a cost such as $||V-V_{max}||^2$ with $V_{max}$ constant would cause the velocity cost to change too much according to the situation, thus making it difficult to tune its weight. A small weight would cause the MPC to be very conservative when it could drive at higher speeds, while a bigger weight could cause safety problems when the velocity should be low, as the optimal solution might be to drive fast rather than to respect some safety constraints.

Therefore, we introduced a local velocity planner that gives a heuristic velocity $V_{heur}$ defined by Equation~(\ref{eq:Vheur}). It takes into account the current speed of the vehicle $V$, a predefined maximum speed allowed $V_{max}$ and the future coordinates $(x_i,y_i)_{i \in \mathbb{N}}$ of the reference path. The strategy consists in increasing the actual speed of the vehicle $V$ by $\Delta V$ while $V_{max}$ is not reached, except if the curvature $\gamma_{max}$ of the path in the next $T_{prev}$ seconds leads to a lateral acceleration $a_y=\gamma_{max} V^2$ over the authorized limit for a kinematic bicycle model of $0.5\mu g$. 
$V_{max}$ is defined specifically for straight lines, to ensure that the vehicle is able to come to a full stop if necessary at the end of the prediction horizon of the MPC. 
%
\begin{eqnarray}
	V_{heur} & = & \min (\sqrt{0.5 \mu g R_{min}}, V_{max}, V + \Delta V)
\label{eq:Vheur}
\end{eqnarray}

\textit{Remark:} the local velocity planner only gives a heuristic velocity for the motion planner, in order to obtain a velocity cost that does not depend to much on the velocity. The goal of the MPC is not to track it perfectly. Therefore, it does not need to be dynamically feasible. 

%
%

\subsection{The Obstacle Manager} 
\label{ssec:obstacle}
This module defines a parabola around each surrounding obstacle in a similar way than \cite{Ziegler2014a}. More precisely, each obstacle $o$ is represented by a parabola $p_0$ with directrix parallel to the reference path at its point closest to $o$; we then choose the parameters of $p_0$ such that it is the minimal parabola containing all vertices of $o$. The interior points of the parabola defines then a region that is forbidden for the vehicle, and will be used as a constraint in the local MPC planner. Please note that only static obstacles are taken into account in this paper. Moreover, the problem of obstacle detection is not considered.

\subsection{The Local MPC Planner} 
\label{ssec:planner}
The local MPC planner updates the reference trajectory every 100ms (10Hz) to take into account changes in the environment, such as obstacles. The concept of Model Predictive Control is to have a \textit{model} of the plant to \textit{predict} the future outputs of the system \cite{Camacho1999}. The MPC solves an optimal control problem which goal is to minimize a cost function $J$ subject to some operational constraints, notably due to vehicle dynamics and obstacle avoidance. Normally, only the first optimal control value computed by the MPC is sent to the plant at time $t$. However, the MPC control outputs here are the acceleration of the vehicle and the steering rate which do not correspond to the actual control inputs which will be used in our simulations, namely the wheel torques and the steering angle (see Section~\ref{sec:SimModel}). Therefore, we used the predicted trajectory computed by the MPC as the reference for the low-level controllers.

At each refresh time, the optimization process is running again with a receding horizon control. The main settings of a MPC is the design of the cost function $J$, the control horizon $N_u$, the prediction horizon $N_y$, and the discretization time step $\Delta t_u$. In our case, the prediction horizon is set to $T_H=3$s, the control discretization time step $\Delta t_u = 0.2$s, and the state prediction horizon $N_y$ is set equal to the control prediction horizon $N_u =T_H/\Delta t_u +1$.

Let $u_1$ and $u_2$ denote respectively the acceleration of the vehicle and the steering rate. Let $U_{t+k,t}=[u_1,u_2,obs_{tol},x_{tol},y_{tol}, \delta_{tol}]$ and $\xi_{t+k,t} = [s,X,Y,V,\psi, \delta]$ denote respectively the control input and the state of our MPC at time $t+k$, predicted at time $t$, with $k=0..N_u$ for $U_{t+k,t}$ and $k=0..N_y$ for $\xi_{t+k,t}$. \\

\subsubsection{Cost function}
We set $e_v=(V_i-V_{heur_i})_{i \in [1:N_y]}$ to be the speed cost; $e_\delta = (\delta_i)_{i \in [1:N_y]}$ and $\dot{e}_{\delta} = (\dot{\delta}_j)_{j \in [1:N_u]}$ respectively the cost on the steering angle and the steering angle variation; $e_{X_{tol}} = (x_{{tol}_j})_{j \in [1:N_u]}$ and $e_{Y_{tol}} = (y_{{tol}_j})_{j \in [1:N_u]}$ the cost of the slack variables associated respectively with the longitudinal and lateral offset to the reference trajectory; $e_{Obs_{tol}} = (obs_{{tol}_j})_{j \in [1:N_u]}$ the cost of the slack variables associated with the obstacle constraints; $e_{\delta_{tol}} =(\delta_{{tol}_j})_{j \in [1:N_u]}$ the cost of the slack variable associated with the validity of the kinematic bicycle model constraint.

We define the following cost function for a control sequence $U=[U_{t,t}, U_{t+1,t},...,U_{t+N_u-1,t}]$ and a state sequence $\xi = [\xi_{t+1,t}, \xi_{t+1,t},...,\xi_{t+N_y,t}]$:
\begin{eqnarray}
J(U, \xi) & = & ||e_v||_{Q_v}^2 + ||e_\delta||_{Q_\Delta}^2 + ||\dot{e}_{\delta}||_{Q_{\dot{\delta}}}^2 + ||e_{X_{tol}}||_{Q_{x}}^2 \\
& & + ||e_{Y_{tol}}||_{Q_{y}}^2  + ||e_{Obs_{tol}}||_{Q_{obs}}^2 + ||e_{\delta_{tol}}||_{Q_{\delta_{tol}}}^2 \nonumber
\end{eqnarray}
where $||x||^2_Q=x^TQx$ and $Q_v=q_v I_{N_y}$, $Q_\Delta=q_\Delta I_{N_y}$, $Q_{\dot{\delta}}=q_{\dot{\delta}} I_{N_u}$, $Q_x=q_x I_{N_u}$, $Q_y=q_y I_{N_u}$, $Q_{obs}=q_{obs} I_{N_u}$, $Q_{\delta_{tol}}=q_{\delta_{tol}} I_{N_u}$. For our simulation, we chose $q_v=4$, $q_{\Delta} = 10$, $q_{\dot{\delta}} = 0.2$, $q_x=5$ $q_y=5$, $q_{obs}=100$ and $q_{\delta_{tol}}=100$.

\subsubsection{Constraints}
First, we constraint the state $\xi_{t,t}$ to be equal to the initial state of the vehicle at time $t$. 
Then, the admissible trajectories of the vehicle are the ones that respect the equations of the kinematic bicycle model described in Equation~(\ref{eq:kin_bic_model}), where we add a constraint on the curvilinear abscissa $s$: $\dot{s} = V$.
The actuator are limited to $u_1 \in [-8;+6]$m/s$^2$ and $u_2 \in [-0.5;+0.5]$ rad/s$^2$. 
The road limits are encoded through a soft constraint: the lateral and longitudinal deviations from the reference path have to remain below a certain tolerance margin $x_{tol}$ and $y_{tol}$.
For each obstacle, we have an inequality constraint which is similar to the one described in \cite{Ziegler2014a}.	


%



At last, to ensure the validity of the kinematic bicycle model, we used the criterion derived in \cite{Polack2017a}, \textit{i.e} that the lateral acceleration $a_y$ ($=V^2/R$ where $R$ is the curvature radius) should be lower than $0.5\mu g$. Using Equation~(\ref{eq:delta_kin}), we obtain the maximum authorized steering angle for a given speed $V$:
\begin{eqnarray}
\delta_{max} (V) = \atan\left(\left(\frac{l_f}{l_r}+1\right)\tan\left(\asin\left(\frac{0.5\mu g l_r}{V^2}\right)\right)\right)
\end{eqnarray}
This is then expressed as a soft constraint by introducing the slack variable $\delta_{tol}$:
\begin{eqnarray}
-\delta_{max} + |\delta| \leq \delta_{tol}
\end{eqnarray}


\subsubsection{Numerical resolution}
The numerical optimization was done using ACADO Toolkit. The working principle of the resolution algorithm can be found in \cite{Houska2011}.\\


\subsection{The Low-level Controllers} 
\label{ssec:controller}

The low-level controllers are responsible for translating the trajectory computed by the local planner into actuator inputs. Many recent research on low-level control have been using MPC with complex vehicle and tire dynamics models for computing the actuators inputs (see for example \cite{Falcone2007a}). 
Although this method enables to anticipate the slip and skid of the vehicle, it is computationally too expensive to be computed at a high frequency. The control is thus in open-loop on a long duration (around 50 to 100ms), which is less robust to modeling errors and might also be less stable.

Therefore, in our control architecture, the low-level control is done in closed-loop at a high frequency (100Hz). Longitudinal and lateral control are treated separately.

For the longitudinal controller, we implemented a simple PID controller \cite{Aastrom1995} 
which takes into account the difference between the reference speed $V^r$ given by the local MPC planner and the actual speed of the vehicle $V$, and its derivatives. $V^r$ is more precisely a linear interpolation of the reference speed output given by the MPC planner. If $e=V-V^r$, then we have the following control law:
\begin{eqnarray}
u(t) & = &  - K_P e(t) - K_D \dot{e}(t) - K_I \int_0^t e(\tau)d\tau
\end{eqnarray}
This control value is then dispatched equally on each of the four wheels of the vehicle in case of braking, and each of the two front wheels in case of acceleration.

For the lateral controller, the steering angle applied $\delta$ is composed of an open-loop part $\delta_{ol}$ and a closed-loop part $\delta_{cl}$. The open-loop steering angle $\delta_{ol}$ is computed by integrating the first value of the control obtained by the MPC planner, which corresponds to the steering rate. This value is only refreshed once every 100ms. The closed-loop steering angle $\delta_{cl}$ is a simple PID control applied on the yaw angle error, projected one MPC time step ahead. 
It is computed every 10ms and uses as reference value for the yaw angle a linear interpolation of the one computed by the MPC planner.


\section{The 9 DoF simulation model}
\label{sec:SimModel}

\begin{figure}[]
	\centering
	\includegraphics[scale=0.4]{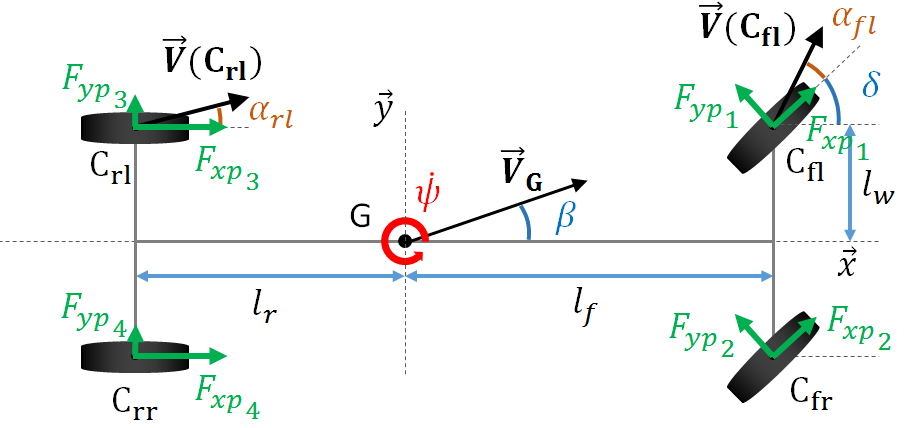}       
	\caption{Simulation model of the vehicle.}
	\label{fig:carSim}
\end{figure}

In this section, we describe the 9 Degrees of Freedom (9~DoF) vehicle model used in order to simulate the dynamic of the vehicle. They correspond to 3~DoF for the whole vehicle ($V_x, V_y, \dot{\psi}$), 2~DoF for the carbody ($\dot{\theta}, \dot{\phi}$) and 4~DoF for the wheels ($\omega_{fl},\omega_{fr},\omega_{rl},\omega_{rr}$). The model takes into account both the coupling of longitudinal and lateral slips and the load transfer between tires. The control inputs of the simulator are the torques $T_{\omega_i}$ applied at each wheel $i$ and the steering angle of the front wheel $\delta$. The low-level dynamics of the engine and brakes are not considered here. The notations are given in Table \ref{tab:notations} and illustrated in Figure~\ref{fig:carSim}.

\textit{Remark: }the subscript $i=1..4$ refers repectively to the front left ($fl$), front right ($fr$), rear left ($rl$) and rear right ($rr$) wheels.

\begin{table}[h]
	\vspace{+0.08in}
	\caption{Notations}
	\label{tab:notations}
	\begin{tabular}{p{1.2cm} p{6.4cm}}
		\hline
		\\
		$\theta$, $\phi$, $\psi$ & Roll, pitch and yaw angles of the carbody \\
		$V_x$, $V_y$ & Longitudinal and lateral speed of the vehicle in its inertial frame \\
		$M_T$ & Total mass of the vehicle\\
		$I_x$, $I_y$, $I_z$ & Inertia of the vehicle around its roll, pitch and yaw axis\\
		$I_{r_i}$ & Inertia of the wheel $i$\\
		$T_{\omega_i}$ & Total torque applied to the wheel $i$\\
		$F_{xp_i}$, $F_{yp_i}$  & Longitudinal and lateral tire forces generated by the road on the wheel $i$ expressed in the tire frame\\
		$F_{x_i}$, $F_{y_i}$ & Longitudinal and lateral tire forces generated by the road on the wheel $i$ expressed in the vehicle frame  $(x,y)$\\
		$F_{z_i}$ & Normal reaction forces on wheel $i$\\
		$l_f$, $l_r$ & Distance between the front (resp. rear) axle and the center of gravity\\
		$l_w$ & Half-track of the vehicle\\ 
		$h$ & Height of the center of gravity\\
		$r_{eff}$ & Effective radius of the wheel\\ 
		$\omega_i$ & Angular velocity of the wheel $i$ \\
		$V_{xp_i}$ & Longitudinal speed of the center of rotation of wheel $i$ expressed in the tire frame\\
		\\
		\hline
	\end{tabular}
\end{table}


Several assumptions were made during the modeling:
\begin{itemize}
	\item Only the front wheels are steerable.
	\item The roll and pitch rotations happen around the center of gravity.
	\item The aerodynamic force is applied at the height of the center of gravity. Therefore, it does not involve any moment on the vehicle.
	\item The slope and road-bank angle of the road are not taken into account.
\end{itemize}

\subsection{Vehicle model}
\label{ssec:CarModel}

The dynamics of the longitudinal, lateral and yaw motions of the whole vehicle are given by the following Equations:
\begin{subequations}
	\begin{eqnarray}
	M_T \dot{V}_x & = & M_T\dot{\psi} V_y + \sum_{i=1}^4 F_{x_i} - F_{aero}\\
	M_T\dot{V}_y & = & - M_T\dot{\psi} V_x +  \sum_{i=1}^4 F_{y_i}\\
	I_z\ddot{\psi} & = & l_f (F_{y_1} +F_{y_2}) - l_r (F_{y_3} + F_{y_4}) \\ \nonumber
	& + & l_w (F_{x_2}+F_{x_4}-F_{x_1}-F_{x_3})
	\end{eqnarray}	
	where the aerodynamic drag forces are $F_{aero} = \frac{1}{2} \rho_{air} C_x S V_x^2$
	with $\rho_{air}$ the mass density of air,  $C_x$ the aerodynamic drag coefficient and $S$ the frontal area of the vehicle. 
	
	The position $(X,Y)$ of the vehicle in the ground frame is given by the following Equations:
	\begin{eqnarray}
	\dot{X} & = & V_x \cos\psi - V_y \sin\psi\\
	\dot{Y} & = & V_x \sin\psi + V_y \cos\psi	
	\end{eqnarray}
	

	\subsection{Carbody model}
	\label{ssec:CarbodyModel}
	The dynamics of the roll and pitch motions of the carbody are given by the following Equations:
	\begin{eqnarray}
	I_x\ddot{\theta} & = & l_w (F_{z_1}+F_{z_3}-F_{z_2}-F_{z_4}) + h \sum_{i=1}^4 F_{y_i}\\
	I_y\ddot{\phi} & = & l_r (F_{z_3} + F_{z_4}) -l_f (F_{z_1} +F_{z_2}) - h \sum_{i=1}^4 F_{x_i}
	\end{eqnarray}	
	where $F_{z_i}$ are damped mass/spring forces depending on the suspension travel $\zeta_i$ due to the roll $\theta$ and pitch $\phi$ angles:
	\begin{eqnarray}
	F_{z_i} = - k_s \zeta_i(\theta, \phi) - d_s \dot{\zeta_i}(\theta, \phi)
	\end{eqnarray}
	
	The parameters $k_s$ and $d_s$ are respectively the stiffness and the damping coefficients of the suspensions.
	\subsection{Wheel dynamics}
	\label{ssec:WheelModel}
	The dynamics of each wheel $i=1..4$ expressed in the tire frame is given by the following Equation:
	\begin{eqnarray} 
	I_r \dot{\omega}_i & = & T_{{\omega}_i}-r_{eff} F_{xp_i}
	\end{eqnarray}
\end{subequations}

\subsection{Tire Dynamics}
\label{ssec:TireModel}
The longitudinal force $F_{xp_i}$ and the lateral force $F_{yp_i}$ applied by the road on each tire $i$ are functions of the longitudinal slip ratio $\tau_{x_i}$, the side-slip angle $\alpha_i$, the normal reaction force $F_{z_i}$ and the road friction coefficient $\mu$: 
\begin{subequations}
	\begin{eqnarray}
	F_{xp_i} & =  & f_x(\tau_{x_i}, \alpha_i, F_{z_i}, \mu)\\
	F_{yp_i} & =  & f_y(\alpha_i, \tau_{x_i}, F_{z_i}, \mu)
	\end{eqnarray}
\end{subequations}

The longitudinal slip ratio of the wheel $i$ is defined as following:
\begin{itemize}
	\item Traction phases ($r_{eff}\omega_i \geq V_{xp_i}$): $\tau_{x_i} = \frac{r_{eff} \omega_i - V_{xp_i}}{r_{eff}|\omega_i|}$
	\item Braking phases ($r_{eff}\omega_i < V_{xp_i}$): $\tau_{x_i} = \frac{r_{eff} \omega_i - V_{xp_i}}{|V_{xp_i}|}$
\end{itemize}
The lateral slip-angle $\alpha_i$ of tire $i$ is the angle between the direction given by the orientation of the wheel and the direction of the velocity of the wheel (see Figure~\ref{fig:carSim}):
\begin{eqnarray}
\small
\alpha_f = \delta - \atan \left(\frac{V_y + l_f \dot{\psi}}{V_x \pm l_w \dot{\psi}}\right) ; \; \alpha_r = - \atan \left (\frac{V_y - l_r \dot{\psi}}{V_x \pm l_w \dot{\psi}}\right)
\end{eqnarray}

In order to model the functions $f_x$ and $f_y$, we used the combined slip tire model presented by Pacejka in \cite{Pacejka2002} which takes into account the interaction between longitudinal and lateral slips on the force generation. Therefore, the friction circle due to the laws of friction (see Equation (\ref{eq:friction_circle})) is respected. Finally, the impact of load transfer between tires is also taken into account through $F_z$. 
\begin{eqnarray}
\label{eq:friction_circle}
||\vec{F}_{xp}+\vec{F}_{yp}|| \leq \mu ||\vec{F}_z||
\end{eqnarray}  

Lastly, the relationships between the tire forces expressed in the vehicle frame $F_x$ and $F_y$ and the ones expressed in the tire frame $F_{xp}$ and $F_{yp}$ are given in Equation~(\ref{eq:frame_F_change}):
\begin{subequations}
	\label{eq:frame_F_change}
	\begin{eqnarray}
	\small{F_{x_i}} & = & \small{(F_{xp_i}\cos\delta_i-F_{yp_i}\sin\delta_i)\cos\phi-F_{z_i}\sin\phi}\\
	\small{F_{y_i}} & = & \small{(F_{xp_i}\cos\delta_i-F_{yp_i}\sin\delta_i)\sin\theta\sin\phi}\\ \nonumber
	& + & \small{(F_{yp_i}\cos\delta_i+F_{xp_i}\sin\delta_i)\cos\theta + F_{z_i}\sin\theta\cos\phi}		
	\end{eqnarray}
\end{subequations}


\section{Results}
\label{sec:results}

\begin{figure}[thpb]
	\centering
	\includegraphics[width=8.5cm, height=5cm]{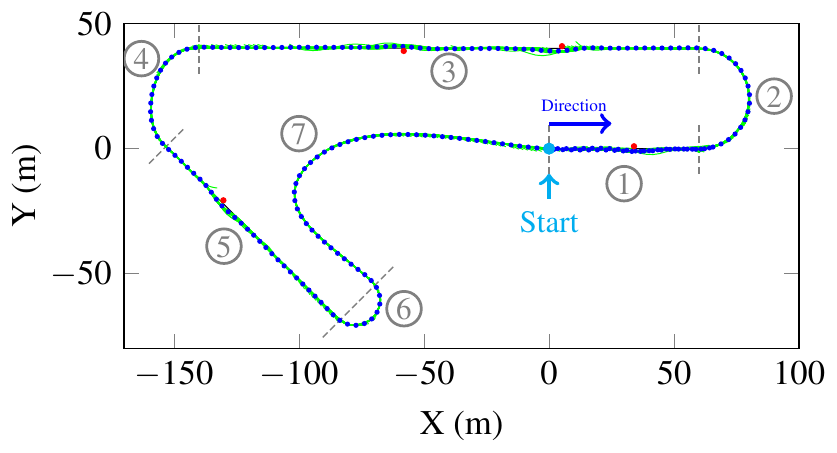} 
	\caption{Trajectories planned by MPC (green lines) and actual trajectory followed by the controller (blue dots), with obstacles (red dots). The numbers indicate different road sections to facilitate the matching with Figures~\ref{fig:Vcompare}, \ref{fig:steering} and \ref{fig:Cmot}.}
	\label{fig:traj_with_obst}
\end{figure}

In order to test the proposed planning and control architecture, simulations were performed in PreScan \cite{Prescanurl} using the 9~DoF vehicle model presented in Section~\ref{sec:SimModel}. PreScan enables to easily generate a traffic scenario in order to test algorithms for autonomous vehicles. The test track was designed to alternate between long straight lines and small radius curvatures (up to 10m only) as shown in Figure~\ref{fig:traj_with_obst}. This is particularly challenging as the constraint on lateral acceleration is rapidly reached in the curves if the vehicle arrives too fast. Therefore, some precautions must be taken in the choice of the parameters of the planning and control architecture. For example, the prediction horizon of the MPC depends on the maximum speed $V_{max}$ as the vehicle needs to be able to decelerate sufficiently in advance to reach the target speed in the curves. The road friction coefficient $\mu$ was chosen equal to 1. The initial speed is $V_{x_0}=15$m/s.

\subsection{Without obstacles}



\begin{figure*}[thpb]
	\vspace{+0.02in}
	\centering
	\includegraphics[width=\textwidth, height=2.7cm]{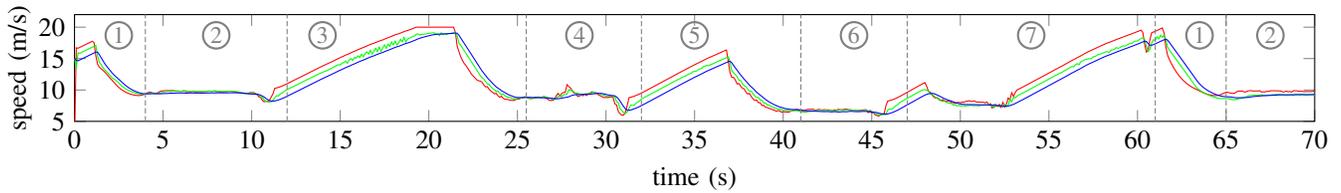}       
	\caption{Comparison between $V_{heur}$ (red), $V^r$ (green) and $V$ (blue), in the case of no obstacles.}
	\label{fig:Vcompare}
\end{figure*}

\begin{figure*}[thpb]
	\centering
	\includegraphics[width=\textwidth, height=4cm]{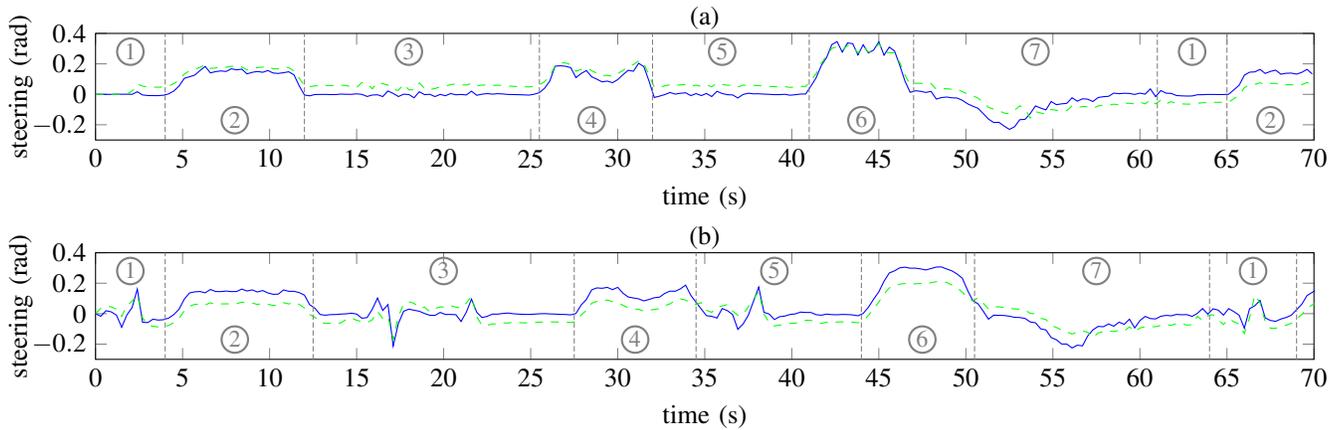}      
	\caption{Total steering angle $\delta$ (blue) and closed-loop steering angle $\delta_{cl}$ (green): (a) without obstacles; (b) with obstacles.}
	\label{fig:steering}
\end{figure*}

\begin{figure*}[thpb]
	\centering
	\includegraphics[width=\textwidth, height=2.8cm]{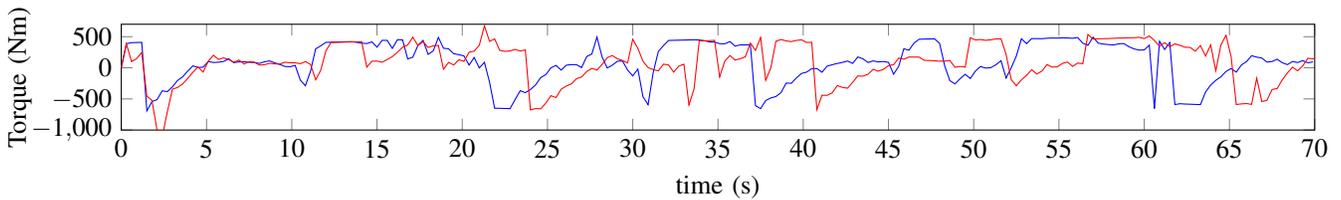}      
	\caption{Wheel torque applied on the front left wheel: without obstacles (blue); with obstacles (red).}
	\label{fig:Cmot}
\end{figure*}

In a first step, we tested the planning and control architecture on the track without obstacles in order to show its capacity to adapt the speed accordingly to the situation. Figure~\ref{fig:Vcompare} compares the heuristic speed $V_{heur}$ computed by the local velocity planner, the target speed $V^r$ computed by the MPC planner and the real speed $V$ of the vehicle. We observe that before entering a curve, the heuristic speed computed by the local velocity planner reduces drastically in order not to exceed $0.5g$ for the lateral acceleration. As there are no obstacles, the speed $V^r$ of the MPC planner is guided by this heuristic, while the longitudinal low-level controller tracks  sufficiently well this reference speed in curves to avoid skipping and skidding. 
The low-level control inputs are shown in Figure~\ref{fig:steering}a and \ref{fig:Cmot} (in blue). 
The absolute value of the lateral error to the reference path given at entry of our system does not exceed $0.4$m and the computational time of the MPC remains always below $100$ms. 
Moreover, we observe in Figure~\ref{fig:steering}a that $\delta_{cl}$ is important due to the slow refresh time of the MPC. 

\pgfplotsset{every axis title/.append style={at={(0.5,0.8)}}}


\subsection{With obstacles}


In a second step, we tested the planning and control architecture on the track with static obstacles in order to show the capacity of the motion planner to cope with planning safe and feasible trajectories. The resulting trajectory is shown in Figure~\ref{fig:traj_with_obst}. The low-level control inputs are shown in Figure~\ref{fig:steering}b and \ref{fig:Cmot} (in red). As expected, the vehicle does almost not slip or skid thanks to the constraint on lateral acceleration. The computational time of the MPC remains below $100$ms. 


\section{Conclusion}
\label{sec:CCL}

We have presented a simple planning and control architecture for an autonomous vehicle under normal driving conditions based on a 10Hz kinematic bicycle MPC and a 100Hz closed-loop PID controller. The main characteristic of our approach is to guarantee the consistency between the two layers by ensuring the validity of the kinematic bicycle model used for planning at any time through a dynamic constraint on the maximal authorized steering angle. Therefore, only feasible trajectories will be generated. This is done at almost no additional computational cost compared to motion planners already using the kinematic bicycle model. The low-level controller enables the architecture to be robust to disturbances and modeling errors. 
Our architecture is also less conservative as it is guided toward the highest speed that keeps the architecture safe thanks to the computation of a heuristic speed based only on the road geometry.

Our planning and control architecture was tested in simulation on a challenging track, with and without obstacles, using a high-fidelity 9~DoF dynamic vehicle model, taking into account both tire friction circles and load transfer. 


Even, if the proposed architecture is able to cope with most normal driving situations, the constraint on lateral acceleration might sometimes be too harsh, especially in emergency situations. Switching to a more complex model at the planning and control level might then be necessary, for example using \cite{AltcheITSC2017}, \cite{Goh2016}. 
This will be in the scope of future work. At last, the impact of low-friction road conditions on the vehicle architecture should also be investigated for further validations.

\bibliographystyle{IEEEtran}
\bibliography{ACC2017}

\end{document}